\documentclass[11pt,leqno]{article}
\usepackage{amssymb,amsmath,theorem,latexsym}
\usepackage[all]{xy}

\theoremstyle{plain} 

\newtheorem{thm}{Theorem}[section]
\newtheorem{theo}{Theorem}

\newtheorem{cor}[theo]{Corollary}
\newtheorem{corr}[thm]{Corollary}
\newtheorem{con}[thm]{Conjecture}
\newtheorem{defi}{Definition}
\newtheorem{prop}[thm]{Proposition}

\newtheorem{claim}[thm]{Claim}
\newtheorem{prop-defi}[theo]{Proposition-Definition}
\newtheorem{lemma-defi}[theo]{Lemma-Definition}
{\theorembodyfont{\rmfamily} \newtheorem{rem}[thm]{Remark}}
{\theorembodyfont{\rmfamily} }
{\theorembodyfont{\rmfamily} }

\newcommand{\op}[1]{\operatorname{#1}}

\newcommand{\fg}[1]{\pi_{1}(#1)}

\newcommand{\Sh}{\boldsymbol{S}\boldsymbol{h}}
\newcommand{\sh}{\boldsymbol{s}\boldsymbol{h}}

\begin{document}

\title {Deformations of large fundamental groups}
\author{B.~de~Oliveira \thanks{Partially supported by NSF
Postdoctoral research fellowship }\and  L.~Katzarkov\thanks{Partially supported by NSF
Career Award and  A.P.Sloan research fellowship } \and 
M.~Ramachandran } 
\date{ }
\maketitle

\bigskip

\tableofcontents

\section{Introduction}

The following notion  of a variety with a large 
fundamental group was introduced by Koll\'ar \cite{Koll97} in his 
study of the Shafarevich uniformization conjecture. 

\begin{defi} A  smooth projective variety $X$ has a 
{\em large fundamental group} if  for every irreducible    
positive dimensional subvariety 
$Z \subset X$  and every resolution of singularities $\op{res}(Z)$  of
$Z$ the image
 $\op{im}[\pi_{1}(\op{res}(Z))
\rightarrow \pi_{1}(X)]$ is an infinite group. 
\end{defi}

\

\noindent
For future reference, it will be convenient to give a special name to
the subvarieties $Z 
\subset X$ which obstruct
the largeness of the fundamental group.

\begin{defi} A subvariety $Z \subset X$ will be called {\em
$\pi_{1}$-small} if for some (and therefore every) resolution of
singularities $\op{res}(Z)$  of 
$Z$ the image 
$\op{im}[\pi_{1}(\op{res}(Z)) \to \pi_{1}(X)]$ is a finite group. 
\end{defi}

\

\noindent
More generally one says that a projective variety $X$ has a {\em
generically large fundamental group} if $X$ contains at most countably
many $\pi_{1}(X)$-small subvarieties.

\

\noindent
In this paper we investigate whether the property of a 
variety having a generically large and large fundamental groups is
stable under   
K\"{a}hler deformations.

The main tools we will use are collected in the following three
theorems.  Their different nature  give us extra 
flexibility to tackle these kind of problems. The statements of the 
theorems are also of independent interest and have 
applications outside this article.

Before we state the theorems we will need to introduce some notation. 
Let $D$ denote a one dimensional complex disk and let $o \in D$ be a
marked point. For a smooth proper K\"{a}hler 
morphism $f : X \to D$ of complex analytic spaces and a point $t \in
D$ denote by $X_{t}$ the fiber $f^{-1}(t)$. Given a finitely
generated group $\Gamma$ we denote by $\mu : \Gamma \to
\widehat{\Gamma}_{\op{uni}}$ the natural homomorphism from $\Gamma$ to
its prounipotent completion. The homomorphism $\mu$ is called the {\em
Malcev representation} of $\Gamma$. The covering of a space $X$
corresponding to the kernel of the Malcev representation of $X$ is
called the {\em Malcev covering of $X$}. A group $\Gamma$ is called
{\em residually unipotent}  if the Malcev representation of $\Gamma$
has finite kernel. 
First we have the following:

\begin{theo} \label{theo-families} Let $f: X \to D$  be a  K\"{a}hler 
deformation of a 
K\"{a}hler manifold $X_{o}$. Suppose that $\pi_{1}(X_{o})$ is
residually unipotent, then there exist a complex analytic space $Y$
which fits in a commutative diagram
\[
\xymatrix@=12pt{X \ar[dr]_-{f} \ar[rr]^-{s} & & Y \ar[dl]^-{p} \\ & D &}
\]
where:
\begin{itemize}
\item Both $p$ and $s$ are proper analytic maps with connected fibers.
\item The fibers of $p$ are normal projective varieties with large
fundamental groups.
\item For each $t \in D$ the fibers of $s_{|X_{t}} : X_{t} \to Y_{t}$
are $\pi_{1}(X_{t})$-small subvarieties of $X_{t}$.  
\end{itemize}
\end{theo}

\

\noindent
The following theorem contains our basic volume estimate and  gives 
a restriction on the  $\pi_1$-small 
subvarieties that occur in a 
K\"ahler deformations of a variety with Stein universal covering.

\begin{theo}  \label{theo-volume}
Let $f: X \to D$ be a  K\"{a}hler deformation of a compact manifold
$X_{o}$ and let $\{ w_{t} \}_{t \in D}$ be a family of K\"{a}hler
metrics along the fibers of $f$. Assume that $X_{o}$ has a  Stein 
universal covering.
Fix any $N\in {\mathbb N}$. Then there exists a smaller disk $o \in
D(N) \subset D$ such that  
\[
\op{vol}_{w_t}(Z)\cdot(\#\op{im}[\pi_1(\op{res}(Z)) \to \pi_1(X_{t})] +1) >
N
\]
for all $t \in D(N)$ and 
every $\pi_{1}$-small subvariety $Z \subset X_{t}$. 

\end{theo}

\

\noindent 
Finally we have the following statement of a purely topological nature

\begin{theo} \label{theo-topological}
Let  $f: X \to T$ be a proper morphism with 
equidimensional fibers between analytic spaces. Then the subset
\[
\{ t \in T | X_{t} \text{ has a $\pi_{1}$-small irreducible
component } \} \subset T
\]
is a closed analytic subset in $T$.
\end{theo}

\

\noindent 
The main objective of the present 
paper is to  use the results stated above to study the 
behavior of large fundamental groups  of K\"{a}hler surfaces under 
deformations. In this direction we obtain the following:

\begin{theo} \label{theo-main} Let $X$ be a K\"{a}hler surface  and
assume that the Malcev covering of $X$ is an infinite covering. 
Then the property of $X$  having a generically large fundamental group
is stable under  small  K\"{a}hler deformations.
\end{theo}

\

\noindent
As a consequence we have the following:

\begin{cor} \label{cor-intermediate}
Let $X$ be a compact  K\"{a}hler surface, let $U$ be a finite
dimensional linear unipotent group, and let $\rho : \pi_{1}(X)
\to U$ be a representation with an infinite image. Suppose that the
covering of $X$ defined by $\ker(\rho)$ 
is Stein. Then any small K\"{a}hler deformation of 
$X$  has a Stein universal covering.
\end{cor}

\noindent
We make the following conjectures:

\begin{con}  Let $X$ be a K\"{a}hler manifold  and let 
$\rho: \pi_{1}(X) 
\to GL(n, {\mathbb C})$ be an infinite linear representation. 
Then having a generically large fundamental group is stable under 
small  K\"{a}hler deformations of $X$.
\end{con}

\begin{con}  \label{con-large}
Let $X$ be a compact  K\"{a}hler  manifold whose
universal covering is Stein.
Then small  K\"{a}hler deformations of $X$ have large fundamental group. 
\end{con}

If the Shafarevich conjecture is correct this implies

\begin{con} \label{con-stein}
If $X$ is a compact K\"{a}hler manifold with a Stein universal cover,
then every small deformation of $X$ also has a Stein universal cover.
\end{con}

\

\noindent
It was  shown in \cite{BK} that  Conjecture~\ref{con-stein} is not
correct without the  
K\"{a}hler hypothesis. We also note that what was lost in that example 
was 
not the property of the fundamental group being large but the 
holomorphic 
convexity.

\

\noindent
The paper is organized as follows. In section two  we describe the
methods and the main
technical tools used in the proofs of the above theorems.
In section three we give the proofs  of  Theorems \ref{theo-volume} and
\ref{theo-main} and some corollaries.

\bigskip

\noindent
{\bf Acknowledgments:} We thank 
F.Bogomolov, M.Gromov, J.Koll\'ar, T.Pantev  and C.Simpson for useful
conversations.

\section{Tools}

In this section we describe the basic techniques used in the proof of
the main theorem. We start with some facts from nonabelian Hodge theory.

\subsection{The relative Shafarevich morphism in families}
\label{subsec-relative}

Since the $\pi_{1}$-small subvarieties are the obstruction to
largeness of the fundamental group, we have to learn how to contract
such subvarieties in the absolute and relative context. In the
absolute case the answer is given by the Shafarevich holomorphic
convexity conjecture which can be stated as follows: 

\begin{con} \label{con-shafarevich} Let $X$ be a smooth complex
projective variety. Then:
\begin{description}
\item[(1)] There exists a  normal analytic space  $\Sh(\widetilde{X})$ 
with  no
compact complex subspaces and a proper morphism  with connected fibers 
$\sh : \widetilde{X} \to \Sh(\widetilde{X})$. 
This morphism is known as the Cartan-Remmert reduction of $\widetilde{X}$.
\item[(2)] $\Sh(\widetilde{X})$ is a Stein space.
\end{description}
\end{con}
\

\noindent
The validity of the Shafarevich conjecture has been recently established
in several classes of examples due to the efforts of F.Bogomolov, 
F.Campana, P.Eyssidieux, L.Katzarkov, J.Koll\'{a}r, T.Napier, R.Narasimhan, 
M.Nori, T.Pantev,  M.Ra\-ma\-chan\-dran (see e.g.  \cite{EY},
\cite{GR},\cite{LR}, \cite{Ka-Pa-Ra}, \cite{Koll97}, \cite{Koll93}, 
\cite{NR}). The strongest result to date is due to Eyssidieux
\cite{EY}, who showed that the Shafarevich conjecture is true for the
covering corresponding to the intersections of the kernels of all
reductive representations of $\pi_{1}(X)$ in some $GL(n,{\mathbb
C})$.

A starting point of most of the recent research on the Shafarevich
conjecture is the important remark of Koll\'ar  \cite{Koll93} and  
\cite{Koll97} who observed  that condition (1)
above is equivalent to:

\begin{description}
\item[$(1')$] There exists a normal variety $\Sh(X)$ and a proper 
morphism  with connected fibers $\sh : X \to \Sh(X)$, which  
contracts exactly $\pi_{1}(X)$-small subvarieties of $X$. 
\end{description}

\

\noindent
Clearly if such a $\Sh(X)$ exists, it will have a large fundamental
group. Therefore property $(1')$ will allow one to reduce
any geometric question related to the fundamental group to a question
about a variety with a large fundamental group. Much work has
been done
on finding various sufficient conditions for the validity of $(1')$.  
The only result in that direction which works in full generality is due
to Campana
\cite{CA1} and Koll\'ar   \cite{Koll97}, \cite{Koll93} who have shown that 
one can always find a proper rational map $\mathsf{sh} : X \dashrightarrow
\mathsf{Sh}(X)$  
with connected fibers which contracts only $\pi_{1}(X)$-small
subvarieties and such that $\mathsf{Sh}(X)$ has a generically large
fundamental group.

\medskip

In what follows we will need Koll\'{a}r's notion of $H$-Shafarevich morphism 
which is a slight generalization of $(1')$ and which we recall next.
We will say that a subgroup 
$\Gamma \subset 
\fg{X}$ is almost contained in $H$ if the intersection $\Gamma \cap H$
has finite index in $\Gamma$ and we will write $\Gamma \lesssim H$. 
The $H$-Shafarevich morphism is  a proper map
with connected fibers $\sh^{H} : X \to  \Sh^{H}(X)$ to a 
normal variety $\Sh^{H}(X)$ 
which contracts exactly the subvarieties  $Z \subset X$ having the 
property that $\op{im}[\pi_{1}(\op{res}(Z))\longrightarrow \pi_{1}(X)]
\lesssim H$.

Recently using ideas and techniques from non-abelian Hodge theory
developed by Corlette \cite{COR}, Gromov and Schoen \cite{GS} and Simpson 
\cite{Simpson}
the following statements were proved (see \cite{Katz}, \cite{Ka-Pa-Ra}).

\begin{thm}[\cite{Katz,Ka-Pa-Ra}] \label{thm-factorization}
Let $X$ be a projective manifold  and let $\rho : \pi_{1}(X)
\to GL(n,{\mathbb C})$  be a linear representation with kernel $H$. Then
there exists a blow-up $\widehat{X}$ of an \'{e}tale finite cover of $X$
and a projective morphism $\op{alb}_{\rho} : \widehat{X} \to
\op{Alb}_{\rho}(X)$ satisfying:
\begin{description}
\item[(i)] $\op{Alb}_{\rho}(X)$ is a smooth fibration of abelian
varieties over a variety of general type.
\item[(ii)] $\op{alb}_{\rho}$ contracts exactly the subvarieties  
$Z \subset \widehat{X}$ having the 
property that 
\[
\op{im}[\pi_{1}(\op{res}(Z))\longrightarrow \pi_{1}(X)]
\lesssim H.
\]
\end{description}
In particular 
the $H$-Shafarevich morphism $\sh^{H} : \widehat{X} \to 
\Sh^{H}(\widehat{X})$ exists and coincides with the Stein
factorization of $\op{alb}_{\rho}$. 
\end{thm}

As an immediate consequence of this theorem one shows

\begin{corr}[\cite{Katz,Ka-Pa-Ra}] \label{cor-shafarevich}
Let $X$ be a K\"{a}hler manifold  and let $\rho : \pi_{1}(X)
\to GL(n, {\mathbb C})$ be an almost faithful representation. Then the 
Shafarevich morphism $\sh : X \to \Sh(X)$ exists and
$\Sh(X)$ and $\mathsf{Sh}(X)$ are birationally isomorphic.
Furthermore if $X$ is a surface, then  the  Shafarevich conjecture 
is true  for $X$. 
\end{corr}

\medskip

\noindent
In order to understand how the varieties with large
fundamental groups vary in families one has to understand how 
the Shafarevich variety $\Sh(X)$  of a projective manifold $X$ changes under
small deformations of the  
complex structure of $X$. We believe that the following Conjecture
should be true:

\begin{con} \label{conj-families} Let $f: X \to D$  be a  K\"{a}hler 
deformation of a  
K\"{a}hler manifold $X_{o}$ over $D$. Let $H \subset \pi_{1}(X) =
\pi_{1}(X_{o})$ be a subgroup for which the $H$-Shafarevich morphism 
$X_{o} \to \Sh^{H}(X_{o})$ exists. Then
there exist 
\begin{description}
\item[a)] a finite base change $\Delta \to D$ branched at $o \in D$, 
\item[b)] a relative Shafarevich variety $p_{H} : 
\Sh^{H}(X\times_{D}\Delta) \to  \Delta$ and a morphism
\[\sh^{H} : 
X\times_{D}\Delta \to
\op{Sh}^{H}(X\times_{D}\Delta),
\]
\end{description} 
so that for all $t
\in \Delta$ the restriction $\sh^{H}_{|X_{t}} : X_{t} \to 
\Sh^{H}(X_{t})$
is the $H$-Shafarevich morphism for $X_{t}$.
\end{con}

Unfortunately, the methods employed in \cite{Katz,Ka-Pa-Ra} for
proving Theorem~\ref{thm-factorization} and
Corollary~\ref{cor-shafarevich} do not extend to the relative
situation. The main problem which prevents one from proving
Conjecture~\ref{conj-families} by the methods of \cite{Katz,Ka-Pa-Ra}
is analytic in nature. More specifically, recall that  for a
projective manifold $X$  the Shafarevich
variety $\op{Sh}^{H}(X)$ corresponding to
the kernel $H = \ker(\rho)$ 
of a reductive  representation $\rho : \pi_{1}(X) \to GL(n,{\mathbb
C})$ is built as a Castelnuovo-de Franchis reduction of $X$ with
respect to a certain system of multivalued holomorphic one forms on $X$
These forms are cnstructed as pullbacks of some constant one forms on
a non-positively curved space by a $\rho$-equivariant  pluriharmonic
map. In the situation when $f : X \to D$ is a smooth projective
morphism it is very hard to ensure the global existence of such
$\rho$-equivariant  harmonic map and even when exists it is hard
to show that this map has reasonable regularity properties. Some
partial results in that direction have been obtained in \cite{KS} but
they require that some very stringent conditions are imposed on the
boundary $\partial X$ of $X$ and do not seem to be immediately
applicable to the situation at hand. In particular, even when one can
prove that the multivalued forms giving rise to $\op{Sh}^{H}(X_{t})$
exist for every $X_{t}$, it seems quite hard to check that these forms
will vary holomorphically with $t \in D$.

On the other hand, there are several instances in which one can prove
Conjecture~\ref{conj-families}  under some mild additional assumption
on the group $H$ or on the K\"{a}hler deformation $f : X \to D$. For
example Theorem~\ref{theo-families} asserts validity of
Conjecture~\ref{conj-families} under the additional condition that $H$
is the kernel of the Malcev representation of $\pi_{1}(X)$. Slightly
more generally we have the following

\begin{claim} \label{claim-nilpotent} Let $f : X \to D$ be a morphism
with smooth projective connected fibers.  Let $U$ be any proalgebraic group
which is a quotient of the prounipotent completion
$\hat{\pi}_{\op{uni}}(X)$ of the fundamental group of $X$. Let $H =
\ker[\pi_{1}(X) \to U]$. Then the Shafarevich morphism $\sh^{H} :
X_{o} \to \Sh^{H}(X_{o})$ exists and Conjecture~\ref{conj-families}
holds for $f : X \to D$ and $H$.
\end{claim}
{\bf Proof.} Let $\Gamma := \op{im}[\pi_{1}(X) \to U]$ and consider
the image $L := \op{im}[h_{1}(X_{o},{\mathbb Z}) \to
\Gamma/[\Gamma,\Gamma]]$. The finitely generated abelian group $L$
corresponds to some abelian variety $A_{o}$ which is a quotient
$\op{Alb}(X_{o}) \twoheadrightarrow A_{o}$ of the Albanese variety of
$X_{o}$. As explained in \cite[Section~2]{Katz-nilpotent} the
existence of a MHS on $\hat{\pi}_{\op{uni}}(X)$ and the strictness
property of  morphism of MHS implies that the Stein
factorization of the natural morphism $X_{o} \to \op{Alb}(X_{o}) \to
A_{o}$ is precisely the $H$-Shafarevich morphism  
\[
\sh^{H} :
X_{o} \to \Sh^{H}(X_{o}).
\] 
On the other hand, since $f : X \to D$ is smooth, projective and with
connected fibers we may assume that (after maybe shrinking $D$) $f$
admits an analytic section and that we have a relative Albanese
variety $\op{Alb}(X/D) \to D$ and a relative Albanese morphism
\[
\xymatrix@=12pt{X \ar[dr]_-{f} \ar[rr]^-{\op{alb}} & &
\op{Alb}(X/D)\ar[dl] \\ & D &} 
\]
Let $A \to D$ be the family of abelian varieties which is the quotient
of $\op{Alb}(X/D)$ corresponding to $L$. Again by
\cite[Section~2]{Katz-nilpotent} the Stein factorization of the
composition $X \to \op{Alb}(X/D) \to A$ restricts to the
$H$-Shafarevich morphism on each $X_{t}$. The claim is proven. \hfill $\Box$

\

\medskip

\noindent
Another special case of Conjecture~\ref{conj-families} can be obtained
as follows:

\begin{claim} \label{claim-global} Let $f : X \to D$ is a smooth
projective morphism with connected fibers and suppose that there exist
a morphism with connected fibers of smooth projective varieties $\pi: M \to
S$ and a fiber product diagram
\[
\xymatrix{
X \, \ar[r] \ar[d]_-{f} & M \ar[d]^-{\pi} \\
D \, \ar@{^{(}->}[r] & S. 
}
\]
Assume further that $H$ is the kernel of the pullback of some
finite dimensional linear representation $\rho : \pi_{1}(M) \to
GL(n,{\mathbb C})$ of $\pi_{1}(M)$ via the
natural map $\pi_{1}(X) \to \pi_{1}(M)$. Then there exists a blow-up
$\widehat{X}$ of a finite \'{e}tale covering of $X$, so that $H$-Shafarevich
morphism exists for $\widehat{X}_{o}$ and Conjecture~\ref{conj-families} holds
for $\hat{f}: \widehat{X} \to D$ and $H$.
\end{claim}
{\bf Proof.} For any blow-up
$\widehat{X}$ of a finite \'{e}tale covering of $X$,
Theorem~\ref{thm-factorization} implies that the $H$-Shafarevich morphism for
$\widehat{X}_{o}$ exist. To check the validity of 
Conjecture~\ref{conj-families} note that by \cite[Theorem~3]{Simpson}
we can deform $\rho$ to a representation $\lambda$ underlying a complex
variation of Hodge structures and so the restriction of $\lambda$ to
each $\pi_{1}(X_{t})$ will be a representation underlying a complex
variation of Hodge structures on $X_{t}$. Now consider the Stein
factorization $s : X \to Y$ of the 
holomorphic horizontal map corresponding to the periods of these
variations of Hodge structures. As argued in \cite{NRL}, the maps
$s_{|X_{t}} : X_{t} \to Y_{t}$ will be precisely the
$\ker(\lambda_{|\pi_{1}(X)})$-Shafarevich morphisms on each
$X_{t}$. By replacing $X$ with a blow-up $\widehat{X}$ of a finite
\'{e}tale covering of 
$X$ we may assume that $Y$ is smooth. Finally, as explained in the proof of
\cite[Theorem~4.4]{Katz}, one has a dihotomy - either
$s_{|\widehat{X}_{t}}$ is 
in fact the $H$-Shafarevich  on $\widehat{X}_{t}$ or $\lambda$ is
defined over a 
number field and some finite direct sum of conjugates of $\lambda$
underlies a polarized integral variation of Hodge structures ${\mathbb
V}$ of weight one on $Y$. Furthermore the strictness property of maps
of MHS implies \cite[Theorem~4.4]{Katz} that the relative Albanese morphism
$a : X \to A$ for the family of abelian varieties $A$ on $Y$
corresponding to ${\mathbb 
V}$ contracts precisely the subvarieties $Z \subset X$ for which 
$\op{im}[\pi_{1}(\op{res}(Z) \to \pi_{1}(X)] \lesssim
\ker(\rho_{|\pi_{1}(X)})$. Thus $\sh^{H} : X \to \Sh^{H}(X)$ exists
and coincides with the Stein factorization of $a$. The claim is
proven. \hfill $\Box$

\

\medskip

\begin{rem} \label{rem-generic} It is worth observing that the
relative version of Koll\'{a}r's Shafarevich map does exist and is
relatively easy to construct. In other words, for any analytic map $f
: X \to D$ whose fibers are connected K\"{a}hler manifolds and any
subgroup $H \subset \pi_{1}(X)$ we can find
a complex analytic space $\mathsf{Sh}(X)$ which fits in a commutative
diagram 
\[
\xymatrix@=12pt{X \ar[dr]_-{f} \ar@{-->}[rr]^-{\mathsf{sh}^{H}} & &
\mathsf{Sh}^{H}(X) 
\ar[dl]^-{p} \\ & D &} 
\] 
so that $\mathsf{sh}^{H}_{|X_{t}} : X_{t} \dashrightarrow
\mathsf{Sh}^{H}(X)_{t}$ is a Koll\'{a}r's Shafarevich map for all $t \in
D$. Indeed, in order to construct $\mathsf{Sh}(X)$ we only have to
repeat the argument in the proof of \cite[Theorem~3.6]{Koll97}. The
only ingredient in this proof which does not automatically carry over
to the analytic situation is the existence of a locally topologically
trivial family of normal cycles on $X$. This however is guaranteed in
the relative situation (possibly after a finite base change) by 
\cite{pourcin}.

Note also that since by construction $\mathsf{sh}^{H}$ fails to
contract at most countably many compact subvarieties $Z \subset X$ for
which $\op{im}[\pi_{1}(\op{res}(Z)) \to \pi_{1}(X)] \lesssim H$, it
follows that for all but countably many $t \in D$ the rational map
$\mathsf{sh}^{H}_{|X_{t}} : X_{t} \dashrightarrow 
\mathsf{Sh}^{H}(X)_{t}$ will actually coincides with the Shafarevich
morphism for $X_{t}$. If in addition $H$ is the kernel of a finite
dimensional linear representation of $\pi_{1}(X)$ we can apply
Theorem~\ref{thm-factorization} to conclude that the $H$-Shafarevich
morphism exists fiber by fiber and so in each of the countably
many fibers $\mathsf{Sh}^{H}(X)_{t}$ of  $\mathsf{Sh}^{H}(X)$ which do
not equal $\Sh^{H}(X_{t})$ we will have at most finitely many
subvarieties $Z \subset \Sh^{H}(X_{t})$ for which
$\op{im}[\pi_{1}(\op{res}(Z)) \to \pi_{1}(\Sh^{H}(X_{t}))] \lesssim H$.

\end{rem}

\subsection{Volume tools}

In this section we study the K\"{a}hler deformations of compact K\"{a}hler 
manifolds
$X_{o}$ whose universal coverings are Stein. A  K\"{a}hler deformation of 
$X_{o}$
means a  K\"{a}hler manifold $X$ and  a smooth proper 
morphism $f : X \to D$ to the one dimensional complex disk $D$ with 
$X_{o} = f^{-1}(o)$, s.t. 
the 
induced K\"{a}hler metric on the fibers $X_{t}$ varies smoothly. 

>From now on $\widetilde{X}$ will denote the universal covering of $X$ and 
$\sigma : \widetilde{X} \to X$ denotes 
the respective
covering map. If $\widetilde{Z}$ is a compact subvariety of
$\widetilde{X}$ we denote by $Z$ its 
reduced image by $\sigma$. We denote the degree of the mapping 
$\sigma_{|\widetilde{Z}} : \widetilde{Z} \to Z$ by $n_{\widetilde{Z}}$. 
Clearly the fact that $\widetilde{Z}$ is compact is equivalent to  
$\op{im} [\pi_1(\op{res}(Z))\to \pi_1(X)]$ being finite. We denote the 
order of this subgroup of 
$\pi_1(X)$ 
by $\#\op{im}[\pi_1(\op{res}(Z)) \to
\pi_1(X)]$. We introduce a notation:

\begin{defi} \label{defi-volume} Let $X$ be a complex manifold and let
$i : Z \hookrightarrow X$ be a compact subvariety of dimension $d$.
If $w$ is a K\"{a}hler form on $X$ define
\[
\op{vol}_w(Z)= \frac {1}{d!}\int_{Z-\op{Sing}Z}(i^*w)^d= \langle [w]^d
, [Z]\rangle,
\]
$[w] \in H^2(X, {\mathbb C})$ is the K\"{a}hler
class and $[Z]\in H_{2d}(X,{\mathbb Z})$ is the homology class of
$Z$. Similarly if $L$ is a positive line bundle on $X$ define
\[
\op{vol}_{L}(Z)= \int_{Z-\op{Sing} Z}dV_L= \frac 
{1}{d!}\int_{Z-\op{Sing} Z}(i^{*}\nu)^d = \deg_Z(L) 
\]
where $dV_L$ is the volume element coming from a metric associated to a 
curvature form $\nu$ for $L$. 
\end{defi}

\

\bigskip

\noindent
Let $f: X \to D$ be a deformation of the compact 
manifold $X_{o}$ over the disc $D$. Then we also have a 
deformation $\tilde{f} : \widetilde{X} \to D$ of the 
universal covering $\widetilde{X}_{o}$ of $X_{o}$. This 
deformation is naturally  $C^\infty$ trivial. Pick a general point 
$x(o)\in \widetilde{X}_{o}$ and let 
$\Sigma_{x(o)} \subset \widetilde{X}$ be a fundamental domain for the
$\pi_{1}(X_{o})$ action on $X_{o}$
containing $x(o)$. 
Using the $C^\infty$ trivialization $\widetilde{X} \cong X_{o}\times D$ of 
$\tilde{f} : \widetilde{X} \to D$ we can extend $x(o)$ to a smooth
section $x : D \to \widetilde{X}$ of the structure map $\tilde{f}$.
Also since $\pi_{1}(X) = \pi_{1}(X_{o})$ the subset $\Sigma :=
\Sigma_{x(o)}\times D \subset \widetilde{X}$ will be a fundamental
domain for the $\pi_{1}(X)$ action on $\widetilde{X}$ containing $x(o)$.
Therefore the subset $\Sigma_{x(t)} = \Sigma \cap X_{t}$ will be a
fundamental domain for the $\pi_{1}(X_{t})$ action on
$\widetilde{X}_{t}$ containing $x(t)$. with this notation established
we are now ready to prove:

\begin{claim} \label{claim-volume}
Let $f:X \to D$ be a K\"{a}hler deformation of 
$X_{o}$ and let $\{ w_{t} \}_{t \in D}$ be a smooth family of
K\"{a}hler metrics along the fibers of $f$.  Assume that $X_{o}$ has 
a Stein universal covering.
Fix any $N \in {\mathbb N}$. Then there exists a smaller disk $D(N)
\subset D$ centered at $o$ so that
\[
\op{vol}_{\tilde{w}_t}(\widetilde{Z}) > N
\]
for all compact subvarieties $\widetilde{Z} \subset \widetilde{X}_t$ 
where $\widetilde{X}_{t}$ is 
the universal covering of $X_t$ and $\tilde{w}_t$ is the pullback of $w_{t}$.
\end{claim}
{\bf Proof.}  The deformation $\tilde{f} : \widetilde{X} \to D$ 
of the Stein manifold $X_{o}$ 
is locally pseudo-trivial \cite{ANVAS}, 
i.e. for any relatively compact  open subset 
$U_{o} \subset \tilde X_{o}$ we can find an open set 
$ U \subset  X$ such that $ U \cap \widetilde{X}_{o}= U_{o}$, and for
which the family $f_{|U} : U \to D$ defines a holomorphically
trivial deformation of $U_{o}$. The open
sets $U_t=  U \cap \widetilde{X}_{t}$ are obtained by flowing $U_{o}$ along a 
holomorphic  
vector field. Choose $U_{o}$ so that it contains a very big ball  
$B_{x(o)}(r)$,  
$r>>0$, around a point $x(o)\in \widetilde{X}_{o}$. 
Since a fundamental domain $\Sigma_{x(o)}$ containing $x(o)$ is relatively 
compact we 
have $\Sigma_{x(o)}
\subset B_{x(o)}(r')$, for some $r'$. Since $U_{t}$, $x(t)$,
$\Sigma_{x(t)}$ and $w_{t}$ vary smoothly with $t$ we can shrink $D$
so that for all $t \in D$ we have
\begin{description}
\item[1)]  $B_{x(t)}(r) \subset U_t$, and 
\item[2)]  $\Sigma_{x(t)} \subset B_{x(t)}(r').$
\end{description}

Condition 1) implies that all compact subvarieties $\widetilde{Z}$ of 
$\widetilde{X}_{t}$ which intersect the
fundamental domain $\Sigma_{x(t)}$ can not be contained in
$B_{x(t)}(r)$. This 
follows 
from $U_{t}$ being
biholomorphic to an open subset of the Stein space
$\widetilde{X}_{o}$. 
Condition 2) states 
that any subvariety
$\widetilde{Z} \subset X_t$ which intersects $\Sigma_{x(t)}$ 
must have $[r-r']$ (the 
integer 
 part of $r-r'$) disjoint pieces of the type $\widetilde{Z} \cap 
B_{x(t)}(1)$.  

Recall the standard fact saying that for any K\"{a}hler manifold $X$ of
bounded geometry, any $k$-dimensional subvariety
$S \subset X$ and any
points $x \in X$ one has
$\op{vol}_{2k}(S\cap B_x(1)) > C > 0$, where 
$C$
depends only on the bound for the geometry of $X$. Since every
subvariety $\widetilde{Z} \subset X_{t}$ can be translated to asubvariety
intersecting $\Sigma_{x(t)}$ the claim follows 
from:
\[
\op{vol}_{\tilde{w}_t}(\widetilde{Z}) > \sum_{i= 1}^{[r-r']}
\op{vol}_{
\tilde{w}_t}(\widetilde{Z} \cap
B_{x(t)}(1)) > [r-r'].C
\]
and the fact that $r$ can be made as big as wanted by making $D$ 
sufficiently small. Notice that
the minimum 
 $C= \op{inf}_{t\in D}\{C_t\}$
of the geometry bounds $C_t$ described above is not zero
 since the metrics vary smoothly.
\ \hfill $\Box$

\subsection{Topological tools}

In this subsection we prove Theorem~\ref{theo-topological}.  
In the proof we use the topological 
stratification and the "lower semicontinuity" property of the homotopy  
groups of the fibers of a  proper morphism.

\bigskip

\noindent
{\bf Proof of Theorem~\ref{theo-topological}.}  Let $f : X \to T$ be a
proper morphism with equidimensional fibers between analytic spaces.
Let $X_{t} = \cup^{k_t}_{l= 1} X_{t}^{l}$
be the irreducible decomposition of the fiber $X_t$. 
Let $o \in T$ be a point such that the fiber $X_{o}$ does not
have any $\pi_{1}$-small irreducible component. We have 
\[
\#\op{im}[\pi_1(\op{res}(X_{t}^{l}) \to \pi_1( X)]= \infty.
\]
for all $t \in U$ and all $1 \leq l \leq k_{t}$. Due to
\cite[Theorems~2.2; 3.3; 4.14]{Verd76} and the fact that $f$ is 
proper we can find an  open 
neighborhood $U$ of $o \in T$ and a finite 
stratification 
that is $U= \cup_i S_i$ of $U$ by
locally closed irreducible 
subsets  so that 
$f_i := f_{|f^{-1}(S_{i})} : f^{-1}(S_i) \to S_i$ are topologically trivial 
fibrations.  Now pick  a neighborhood $o \in V\subset U$  
s.t. $S_i \cap V \neq \varnothing$ iff $o  
\in \overline S_i$. Consider the morphisms $\bar{f}_{i} :=
f_{|f^{-1}(\overline{S}_{i})} : f^{-1}(\overline{S}_{i}) \to
\overline{S}_{i}$. 
The theorem will follow if we show that the general 
fiber component of each of the 
morphisms  $\bar{f}_{i}$ is not $\pi_{1}$-small.

Pick an irreducible component of $f^{-1}(\overline{S}_{i})$
that dominates $\overline{S}_i$. Let
$\widehat{f^{-1}(\overline{S}_{i})}$ and $\widehat{\overline{S}_{i}}$
denote the normalizations of $f^{-1}(\overline{S}_{i})$ and
$\overline{S}_i$ 
respectively. We get an induced  morphism of normal irreducible 
varieties  $\hat{\bar{f}_{i}} : \widehat{f^{-1}(\overline{S}_{i})} \to 
\widehat{\overline{S}_{i}}$. Notice that the 
irreducible components of the fiber  $\hat{\bar{f}_{i}}^{-1}(t)$  are 
partial  normalizations of some irreducible components of the fiber 
$X_{t}= f^{-1}(t)$.
 Hence the hypothesis of the theorem  is still valid for the morphism 
$\hat{\bar{f}_{i}}$. By construction of the stratification we know that
if $S_{j} \neq S_{i}$ are two strata such that  $S_{j} \subset 
\overline{S}_{i}$ is of minimal codimension then the fiber of 
$\hat{\bar{f}_{i}}$ over a point in $S_{j}$ will be a deformation retract 
of the fiber
of $\hat{\bar{f}_{i}}$ over a point in $S_{i}$. In particular if $Y$ is
the normalization of a component of $\hat{\bar{f}_{i}}^{-1}(t)$, $t \in
S_{i}$, then $\op{im}[\pi_{1}(Y) \to \pi_{1}(X)]$ coincides with 
$\op{im}[\pi_{1}(Z) \to \pi_{1}(X)]$, where $Z$ is the normalization of 
some component of $\hat{\bar{f}_{i}}^{-1}(s)$ for some point $s \in
S_{j}$. Due to this we may assume that $o \in S_{i}$ and so
it suffices to check the validity of the theorem for a
$f: X \to T$ satisfying: (a) $f$ is a  proper morphism of irreducible 
normal complex 
spaces which; (b) $f$ a topological fiber bundle over a dense open 
subset $\dot{T}$ of $T$; (c)  $o \in \dot{T}$. But in this situation the
statement of the theorem is obvious since $f$ is topologically locally
trivial around the point $o$. \hfill $\Box$

\bigskip

\noindent
For future reference we also give a proof of the following topological result 
which is  a 
straightforward generalization  of the 
results of \cite{KoMiMo91}, \cite{CA}, and  \cite{Koll97}.

\

\begin{prop} Let $X$ be a normal projective variety. Then $\pi_1(X)$ 
is finite if and only if $X$ is chain connected via $\pi_1$-small 
curves. 
\end{prop}
{\bf Proof.} One implication is obvious. Suppose that $X$ is chain 
connected 
via $\pi_1$-small
 curves. By \cite[Theorem~1.8]{Koll97}  there exists an open subset 
$\dot{X}$ and a 
proper  morphism
$f: \dot{X} \to Z$, the Shafarevich map, such that if $z \in Z$ then 
$X_z= f^{-1}(z)$ 
is $\pi_1$-small. This map has also the property that if $z \in Z$ is 
very 
general all normal cycles which are $\pi_1$-small  and intersect $X_z$ 
must be contained in $X_z$. Pick a point $x \in X_z$, with $z$ very 
general. Let  $y$ be any other point in $X$. There is a chain of 
$\pi_1$-small  curves $\{C_1,C_2, \ldots ,C_k\}$ with $C_i \cap C_{i+1} \neq 
\varnothing$, $x \in C_1$ and $y \in C_k$. 
The property of the Shafarevich map mentioned above implies that $C_1 
\subset 
X_z$, which in turn implies $C_2 \subset X_z$ and so on and hence $y 
\in X_z$. This implies $f(X)= point$ and we are done. \hfill $\Box$

\section{Theorems on generically large fundamental groups}

\subsection{Proof of the main theorem}

In this subsection we prove Theorem~\ref{theo-volume} and
Theorem~\ref{theo-main}

\

\noindent
{\bf Proof of  Theorem~\ref{theo-volume}.}
Let $X_{o}$ be a compact K\"{a}hler manifold
 with a Stein universal cover and let
$f: X \to D$ be a  K\"{a}hler deformation of $X_{o}$. Let as before $\{ w_{t}
\}_{t \in D}$ denote a smooth family of K\"{a}hler metrics along the
fibers of $f$. Fix $N \in {\mathbb N}$. According to
Claim~\ref{claim-volume} we can find a smaller disk $o \in D(N) \subset
D$ so that  for any $t \in D(N)$ and all compact subvarieties 
$\widetilde{Z}\subset \widetilde{X}_{t}$  we
have $\op{vol}_{\tilde{w}_t}(\widetilde{Z}) > N$.
The theorem now follows from the obvious identities
\[
\op{vol}_{\tilde{w}_{t}}(\widetilde{Z})=
n_{\widetilde{Z}}\cdot\op{vol}_{w_{t}}(Z)
\]
and
\[
n_{\widetilde{Z}}= \#\op{im}[\pi_1(\op{res}(Z)) \to \pi_1(X_t)].
\]
\

\ \hfill $\Box$

We also have some immediate corollaries:

\begin{corr}  \label{cor-Lbound}
Let $X_{o}$ be a compact K\"{a}hler manifold
with a Stein universal
cover. Suppose
$f: X \to D$ is a deformation of $X_{o}$ and let $L \to X$ be a
$f$-ample line bundle. Then for any $N\in {\mathbb N}$ there exists a
smaller disk $o \in D(N) \subset
D$ so that  for any $t \in D(N)$ and all compact subvarieties 
$\widetilde{Z}\subset \widetilde{X}_{t}$  we
have
\[
\deg_{L_t}(Z)\cdot (\#\op{im}[\pi_1(\op{res}(Z)) \to \pi_1(X_t)]+1) > N 
\]
where $L_t= L_{|X_t}$.
\end{corr}
{\bf Proof.} It is straightforward consequence of 
Theorem~\ref{theo-volume} applied to the family of K\"{a}hler metrics
corresponding to the curvature form of relative ample line bundle $L$.
\hfill $\Box$

\begin{corr}   \label{cor-bound}
Let $X_{o}$ be a compact K\"{a}hler manifold 
with a Stein universal cover and let
$f: X \to D$ be a  K\"{a}hler deformation of $X_{o}$. Suppose 
that the order of all finite subgroups in $\pi_{1}(X_{o})$ is bounded
from above. Then for any $N \in {\mathbb N}$ we can find a smaller 
disk $o \in D(N) \subset
D$ so that  for any $t \in D(N)$ and all $\pi_{1}$-small subvarieties $Z
\subset X_{t}$ we have $\op{vol}_{w_{t}}(Z) > N$.
\end{corr}
{\bf Proof.}  Let $b$ be the upper bound for the order of the finite 
subgroups of 
$\pi_1(X_{o})$. By Theorem~\ref{theo-volume} we can find a disk 
$o \in D' \subset D$ so that for any $t \in D'$ and any 
$\pi_{1}$-small subvariety $Z \subset X_{t}$ we have
\[
\op{vol}_{w_t}(Z)\cdot(\#\op{im}[\pi_1(\op{res}(Z)) \to  \pi_1(X)] +1) >
N\cdot b.
\] 
But  by assumption  $\# \op{im}[\pi_1(\op{res}(Z))\to \pi_1(X)] < b$ and
so the corollary is proven. \hfill $\Box$

\medskip

\noindent
The above corollary explains that for a small K\"{a}hler  (projective) 
deformations 
of a K\"ahler (projective) manifold as in the corollary
the property of having a large fundamental group  can be only broken 
by subvarieties with very big volume (degree). Hence if we
have  an upper bound for the degree of the subvarieties that can violate 
the
property of  a large fundamental group, 
we can show that  large fundamental groups are preserved under small 
K\"{a}hler  
deformations. We use this idea to give a proof of
Theorem~\ref{theo-main}. In fact we can proof the following general

\begin{prop} \label{prop-main} Let $f : X \to D$ be a K\"{a}hler
deformation of a K\"{a}hler surface $X$ and assume that $H \subset
\pi_{1}(X)$ is a subgroup of infinite index  such that the
$H$-Shafarevich morphism for $X_{o}$ exists and
Conjecture~\ref{conj-families} 
holds for $f : X \to D$ and $H$. Assume also that $X_{o}$ has a
generically large fundamental group. Then $X_{t}$ has a generically
large fundamental group for all $t$ in a small analytic neighborhood
of $o \in D$.
\end{prop}

\begin{rem} Note that the previous proposition in combination with
Claim~\ref{claim-nilpotent} implies  Theorem~\ref{theo-main}. Moreover
Claims~\ref{claim-nilpotent} and \ref{claim-global} give other
situations to which we can apply Proposition~\ref{prop-main} to
conclude that the property of having a generically
large fundamental group will be preserved under deformations
\end{rem}

We will
give two different proofs Proposition~\ref{prop-main}. (All surfaces
are assumed to be 
compact and connected unless stated otherwise.)

\

\medskip

\noindent
{\bf First Proof.} The 
existence 
of a relative Shafarevich morphism plus the "lower semicontinuity" of 
the the homotopy
groups in a family gives the first proof.

Let $X_{o}$ be a compact K\"{a}hler surface and let $\rho :
\pi_{1}(X_{o}) \to GL(n,{\mathbb C})$ be an infinite linear
representation with kernel $H$.
Let $f: X \to D$ be a deformation of $X_{o}$. By
Theorem~\ref{theo-families} there exists a relative
$H$-Shafarevich morphism for the family 
$f: X \to D$ possibly after a finite base change. Since the stability of
$\pi_{1}$-small curves cannot be affected by a base change we may assume
without a loss of generality that we have a relative $H$-Shafarevich
morphism 
$\sh^{H} : X \to \sh^{H}(X)$ for $X$ itself. Let $g$ denote the 
product 
morphism 
$g :=  \sh^{H}\times f: X \to \Sh^{H}(X)\times D$.

Since the $\pi_1$-small curves must lie in the fibers of $g$, we need to 
show 
that there exists an open neighborhood $U$, $o \in U \subset D$ s.t. for 
all $t \in 
 U$ the fibers $X_{y,t}= g^{-1}(y,t)$ satisfy:
\[
\op{im}[\pi_1(\op{res}(X_{(y,t)})) \to \pi_1(X)]= \infty.
\]
There are two cases:

\medskip

\noindent
{\bf Case 1.}  All fibers $X_{(y,o)}$, $y \in \Sh^{H}(X)$ are zero
dimensional. The dimension of the 
fibers of 
$f$ is an upper-semicontinuous function in the analytic topology 
of $g(X)$. This implies that there exists an open neighborhood $U$ 
of $o \in \Delta$ such that all fibers for $X_{(y,t)}$  for $t \in U$
are  zero dimensional. Hence there are no $\pi_1$-small curves.

\medskip

\noindent
{\bf Case 2.}  Some $X_{(y,o)}$ are 1-dimensional. From the 
upper-semicontinuouty mentioned above and the fact that $g(X_t)$ can not 
be a point it follows 
that the set $S := \{s\in g(X)| \dim g^{-1}(s)= 1\}$ 
is a  closed analytic subspace of $g(X)$. By 
shrinking  $D$ 
we can assume that all irreducible components $S_i$ of 
$S$ pass through $g
(X_{o})$. This implies that all  possible $\pi_1$-small curves on $X_t$, 
$t\in D$, lie in a family of curves in $ 
X$ containing a curve in $X_{o}$. 
Fix a curve $X_{(y_{o},o)}$ and an irreducible component $S_i$ of $S$. The 
family of curves $g : g^{-1}(S_i) \to S_i$,
is a family whose central fiber contains no $\pi_1$-small curves.
Theorem~\ref{theo-topological} then implies that there is an  
open neighborhood  $y_{o} \in U \subset \Sh^{H}(X)\times D$, 
s.t. over $U\cap S_i$ does not have any not $\pi_1$-small curves. Since 
$S$ has only finitely many irreducible components, 
there is only a finite number of families to consider and so we can
choose $U$ which does not contain any 
$\pi_1$-small curves.  Now since $g(X_{o})$ is compact the result follows. 
\hfill $\Box$

\

\medskip

\noindent
{\bf Second  Proof.} We now give a second proof of
Proposition~\ref{prop-main} but under some additional hypothesis -  we
will assume that that $\pi_{1}(X_{o})$ has finite torsion.
In the second  proof we use  the volume tools (Theorem~\ref{theo-volume}).

\begin{claim} \label{claim-main} 
Let $X_{o}$ be a K\"{a}hler surface  with 
$\rho: \pi_{1}(X_{o}) \to GL(n,{\mathbb C})$
an infinite linear representation. Assume that $X_{o}$ has a Stein
universal cover and that $\pi_{1}(X_{o})$ has finite torsion.
Then any small K\"{a}hler deformation of $X_{o}$ has a  
large  fundamental group.
\end{claim} 
{\bf Proof.}   We  will use the relative 
$\rho$-Shafarevich morphism to control the volume of the $\pi_1$-small 
curves and then use Corollary~\ref{cor-bound} to get a contradiction.

Let $g =  \sh^{H}\times f: X \to \Sh^{H}(X)\times D$ be the same 
morphism as in the previous proof. Since $g$ is proper we can assume 
that by shrinking $\Delta$ the set of points $y \in g(X)$ with 
1-dimensional 
$g^{-1}(y)$  is a finite
union of irreducible closed sets, $\cup_i F_i$. Let 
$X_{F_i}= \op{red}(g^{-1}(F_i))$ be the reduced analytic space associated with 
$g^{-1}(F_i)$. The spaces
$X_{F_i}$ are finite unions    
of irreducible components, hence to simplify the proof we can assume 
$X_{F_i}$ 
is irreducible. Let $\widehat{F}_i$ be the normalization of $F_i$ and 
$g_{\widehat{F}_i}: X_{\widehat{F}_i} \to  \widehat{F}_i$
be the pullback to $\widehat{F}_i$ of $g : X_{F_i} \to {F_i}$.

According to \cite[Theorem~3]{BAR78} the volume of the fibers of 
$g_{\widehat{F}_i}: X_{\widehat{F}_i} \to  \widehat{F}_i$ relative to 
any hermitian
metric over $ X_{\overline{F}_i}$  is locally bounded over $\overline{F}_i$.  
We 
consider the metric induced from the  K\"{a}hler metric on  $ X$. Due to
the fact 
that there are only a finite number of $F_i$  implies that there is 
neighborhood $o \in U \subset D$
s.t all surfaces $X_t$, $t\in U$ have $\pi_1$-small curves with bounded 
volume. Recall that Corollary~\ref{cor-bound}
says that the volume of the $\pi_1$-small curves is larger than any given 
number by making the deformation 
sufficiently small. Hence the result follows by contradiction. \hfill
$\Box$

\bigskip

\noindent

\

Using the previous results we can now prove the stability of the 
Steiness of universal covers mentioned in the introduction:

\bigskip

\noindent
{\bf Proof of Corollary~\ref{cor-intermediate}.} 
This result follows from looking at an intermediate  cover
of $X$, more precisely the infinite cover $X^{\rho}$ defined by 
$\ker(\rho)$. The starting point  are the result of \cite{Katz} stating 
that $X^{\rho}$ is holomorphic convex
and the fact that this follows just from $\rho$ being an infinite linear 
representation. Hence  $X_t^{\rho}$ being holomorphically convex is 
stable under deformations of $X$. 
To get the result one needs only to show that $X_t^{\rho}$
remains Stein for small K\"{a}hlerian deformation of $X$, since 
unramified coverings (even infinite ones)  of Stein varieties are still Stein.

By the previous paragraph, it only remains to show that having a 
$\rho$-large fundamental groups is stable for small K\"{a}hlerian 
deformations of $X$. This however follows from Theorem~\ref{theo-main}
\hfill $\Box$

\

\subsection {Deformations of surfaces and $\pi_1$-small curves}

\

This subsection gives some results on $\pi_1$ -small curves  on 
surface  deformations. We are going to use geometric bounds for the 
degree of 
curves on surfaces plus Theorem~\ref{theo-volume}
 to get  information on the small 
curves that can appear along a deformation. 

We also give a proof for the fact that a small deformation of a surface 
of general type whose 
universal covering is Stein has no rational curves $C$ with $C^2\leq 0$. 
We then 
note that this result plus a consequence of the Shafarevich Conjecture 
would rule out the appearance rational curves.

\

\begin{prop} \label{prop-max}
Let $X_{o}$ be a surface of general type whose universal 
covering is Stein.
Suppose $f: X \to D$ is a deformation of $X_{o}$ over the 
disc. Then for any $N\in {\mathbb N}$ there is a sufficiently small disc 
$o \in D(N) \subset D$ s.t. all 
$\pi_1$-small curves $C \subset X_t$, $t\in D(N)$, satisfy:
\begin{equation}
\op{max} \{ g(C),  C^2,  \#\op{im}[\pi_1(\widehat{C}) \to \pi_1(X_{t})] 
\} > N, \label{eq1}
\end{equation}
where $\widehat{C}$ is the normalization of $C$.
\end{prop}
{\bf Proof.} First note that the canonical class 
$K_{X_{o}}$ is ample since $\pi_1(X_{o})$ 
being large 
rules out the existence of $(-2)$-rational curves. By shrinking the base 
$D$, if necessary, we can choose the relative
canonical bundle $\omega_{X/D}$ as the relative ample line bundle 
required by Corollary~\ref{cor-Lbound}. Thus we can find a disk $o \in D(N)
\subset D$  s.t. all
 small curves $C \subset X_t$, $t\in D(N)$, satisfy:
\begin{equation} \label{eq2}
\max \{  K_{X_t}\cdot C,  \#\op{im}[\pi_1(\widehat{C}) \to \pi_1(X_t)]  
\} > N.
\end{equation} 
Using the results of \cite{Miy-Lu95} in the
spirit of the results of \cite{Bogo77} we know that the degree of a 
curve $C$ 
having  geometric genus
$g(C)$ and selfintersection $C^2$ in a surface of general type $X$  is 
bounded 
by:
\begin{equation} \label{eq3}
K_X\cdot C < 6c_2(X)- 2c_1^2(X)+6g(C)-6 +C^2 
\end{equation}
Since $6c_2(X)-2c_1^2(X)$ is a topological invariant of the surface it 
remains 
constant under deformations. Therefore the combination of (\ref{eq2}) and
(\ref{eq3})  yields the 
desired conclusion. \hfill $\Box$

\bigskip

\noindent
The previous proposition implies that if any two of the three entries 
in (\ref{eq1}) can be bounded the remaining entry must be very high.

\bigskip

\noindent
The study of rational curves $C$ along small deformations is more 
attainable 
because in this case $\#\op{im}[\pi_1(\widehat{C})
\to \pi_1(X_{t})]= 0$ (see the next remark) and so it is not necessary
to give 
a bound for the order of the torsion  subgroups of $\pi_1(X_{o})$.

\begin{rem}  Let $C$ be a rational curve on a complex manifold $X$  and 
let $C'$ be an irreducible component of the preimage of $C$ in the 
universal covering $\sigma :\widetilde{X}\to X$.  Then $C'$ is compact and the 
covering map $\sigma_{|C'}:C' \to C$ is of degree one. To see
this we just notice that there is a map  $\hat{\sigma}:\widehat{C}'\to
\widehat{C}$, between the normalizations of 
$C'$ and 
$C$, that is by construction an unramified 
covering. Since $\widehat{C} \cong {\mathbb P}^{1}$ we have
$\widehat{C}'= {\mathbb P}^1$ and $\deg(\hat{\sigma})= 1$.
\end{rem}

\

\begin{prop} \label{prop-rational}
  Let $X_{o}$ be a surface of general type whose
universal covering is Stein. Then any small deformation of $X_{o}$ has no 
rational 
curves $C$ with $C^2 \leq 0$.
\end{prop}
{\bf Proof.} It follows from the previous remark  that 
$\#\op{im}[\pi_1(\widehat{C}) \to 
\pi_1(X_{t})]= 0$. The
hypothesis of the proposition bounds $C^2$, hence the conclusion is 
obtained 
using Proposition~\ref{prop-max}. \hfill $\Box$

\begin{rem} The Shafarevich Conjecture combined with
Proposition~\ref{prop-rational}
imply that 
if $X$ is a surface with large fundamental group then any small 
K\"{a}hler deformation of $X$ has no rational curves. If the universal 
covering $\widetilde{Y}$ of a surface
$Y$ is holomorphically convex then there is Remmert reduction, 
$r: \widetilde{Y} 
\to S$.  The compact curves $\widetilde{C}$ of $\widetilde{Y}$ are contained in 
the fibers of the reduction map $r$ and hence a $\pi_1$-small     curve 
$C$ of $Y$ must have $C^2\leq 0$. 
\end{rem}

\bigskip

\noindent
B. de Oliveira, University of Pennsylvania and Harvard University. \\ bdeolive@math.upenn.edu, bdeolive@math.harvard.edu

\bigskip

\noindent
L. Katzarkov, UC Irvine. \\ lkatzark@math.uci.edu

\bigskip

\noindent
M. Ramachandran, SUNY Buffalo \\ ramac-m@math.buffalo.edu

\end{document}